\documentclass[11pt]{article}
\usepackage{graphicx}
\usepackage{textcomp}
\usepackage{amsmath,amsfonts,amscd,amsthm,amssymb}
\usepackage{epsfig}
\usepackage{makeidx}

\makeindex
\newcommand{\ncom}{\newcommand}
\newcommand{\be}{\begin{eqnarray*}}
\newcommand{\ee}{\end{eqnarray*}}
\newcommand{\ben}{\begin{eqnarray}}
\newcommand{\een}{\end{eqnarray}}

\newcommand{\sra}{\rightarrow}

\newtheorem{tdf}{Theorem}[section]

\ncom{\ul}{\underline}
\ncom{\beq}{\begin{equation}}
\ncom{\eeq}{\end{equation}}
\ncom{\bea}{\begin{eqnarray*}}
\ncom{\eea}{\end{eqnarray*}}
\ncom{\beqa}{\begin{eqnarray}}
\ncom{\eeqa}{\end{eqnarray}}
\ncom{\nno}{\nonumber}
\ncom{\non}{\nonumber}
\ncom{\ds}{\displaystyle}
\ncom{\half}{\frac{1}{2}}
\ncom{\mbx}{\makebox{.25cm}}
\ncom{\hs}{\mbox{\hspace{.25cm}}}
\ncom{\rar}{\rightarrow}
\ncom{\Rar}{\Rightarrow}
\ncom{\noin}{\noindent}
\ncom{\sz}{\scriptsize}
\ncom{\Sgm}{\Sigma}
\ncom{\psgm}{\sigma^{\prime}}
\ncom{\dt}{\delta}
\ncom{\Dt}{\Delta}
\ncom{\lmd}{\lambda}
\ncom{\Lmd}{\Lambda}
\ncom{\Th}{\Theta}
\ncom{\eps}{\epsilon}
\ncom{\pcc}{\stackrel{P}{>}}
\ncom{\lp}{\stackrel{L_{p}}{>}}
\ncom{\sspan}{{\rm\,span}}
\ncom{\re}{{\rm Re\,}}
\ncom{\im}{{\rm Im\,}}
\ncom{\sgn}{{\rm sgn\,}}
\ncom{\ba}{\begin{array}}
\ncom{\ea}{\end{array}}
\ncom{\integ}[4]{\int_{#1}^{#2}\,{#3}\,d{#4}}
\ncom{\vspan}[1]{{{\rm\,span}\{ #1 \}}}
\ncom{\dm}[1]{ {\displaystyle{#1} } }
\ncom{\ri}[1]{{#1} \index{#1}}

\newtheoremstyle
    {remarkstyle}
    {}
    {11pt}
    {}
    {}
    {\bfseries}
    {:}
    {     }
    {\thmname{#1} \thmnumber{#2} }

\theoremstyle{remarkstyle}

\def \R{{{\rm I{\!}\rm R}}}
\def \D{{{\rm I{\!}\rm D}}}

\begin{document}
\title
  {A  New Modified Newton Method use of Haar wavelet for solving Nonlinear equations }
  \author
{Bijaya Mishra{\footnote {Department of mathematics,Gandhi Institute for Technological Advancement,
Bhubaneswar-752054,India.Email:bijayamishra.math@gmail.com}},
Ambit Kumar Pany {\footnote{Center for Applied Mathematics, Siksha O Anusandhan University, Bhubaneswar, Odisha,
India, 751030. Email:ambit.pany@gmail.com}},
Salila Dutta {\footnote{Department of  Mathematics, Utkal University, Bhubaneswar, Odisha,
India, 751004. Email:saliladutta516@gmail.com}}}
\maketitle 
\begin{abstract}
In this paper, we present a new modified Newton method a use of Haar wavelet formula for solving non-linear equations. 
This new method  do not require the use of the second-order derivative. It is shown that the new method has third-order of 
convergent. Furthermore, some numerical experiments are conducted which confirm our theoretical findings.
\end{abstract}

\setcounter{section}{0}
\section{\normalsize\bf Introduction}
In numerical analysis, finding a solution of  non-linear equation is one of the most attractive problem. In this paper, we
emphasize on an iterative method to find a simple root $\alpha$ of a non-linear equation $f(x)=0$, i.e., $f(\alpha)=0$ and $f'(\alpha)\neq 0$.
Here, we less concern about multiple roots. Newton’s method \cite{OS},\cite{conte} is the well known algorithm
to solve nonlinear equation. It is given by
\begin{equation}
\label{nm}
 x_{n+1}=x_n-\frac{f(x_n)}{f'(x_n)}\,\,\,\,\, n= 0, 1, 2, \cdots
\end{equation}
and it converges quadratically.\\
Earlier, \cite{WF}-\cite{OZ} and \cite{FS1}-\cite{KLW} derived third-order convergence methods based on integral interpretation of Newton's method. Where,
Newton's method derived from different quadrature formulas for the indefinite integral arising from Newton's theorem \cite{DS}
\begin{equation}
\label{nt}
 f(x)= f(x_n)+\int_{x_n}^x f'(t) dt.
\end{equation}
Weerakoon et al.\cite{WF} have approximated the integral part of (\ref{nt}) by trapezoidal rules and derived a variant of Newton's method. It is, further, shown that
this method converges cubically. Subsequently, Frontini et al. \cite{FS} have proposed a third order convergent method 
by approximate the integral by the midpoint rule.
In \cite{HOM}, Homeier has developed a cubically convergent iteration scheme by considering Newton’s theorem for the inverse function. Further,  in \cite{HOM1}, \cite{HOM2}
modified Newton methods are derived for multivariate case. Kou et al. in \cite{KLW} have applied a new interval of integration on Newton’s theorem and arrived
a third-order convergent iterative scheme.\\

Recently, Islam et al.\cite{IAH} have applied Haar wavelet function to derived quadrature rules for indefinite integration. 
In this paper, we modified Newton's theorem by using the quadrature rule proposed by Islam in \cite{IAH}. It is shown that the new method
has third order convergent. Further, the new method did not evaluate second derivative of $f$. The efficiency of the new method
is demonstrated by numerical examples. 

This paper is organized as follows. In Section 2, we discuss a modified Newton's method.
Section 3, we establish convergence analysis for the new method.
finally in section 4, various numerical experiments conduct to confirm our theoretical finding.
\section{\normalsize\bf A Modified Newton's Method}
\setcounter{equation}{0}
To derive the new method, we consider  Newton’s theorem
\begin{eqnarray}
 \label{nt1}
 f(x)= f(x_n)+\int_{x_n}^x f'(\tau) d\tau 
\end{eqnarray}
We use the Haar wavelet function to approximate the integral term of (\ref{nt1}) as
\begin{equation}
 \label{hw}
 \int_{x_n}^x f'(\tau) d\tau = \frac{(x-x_n)}{2M}\sum_{k=1}^{2M} f'\big(x_n+\frac{(x-x_n)(k-0.5)}{2M}\big).
\end{equation}
where $M=2^{J_1}$ and $J_1$ is the maximum level of resolution of Haar wavelets, see \cite{IAH}.\\
Substitute (\ref{hw}) in (\ref{nt1}) to obtain 
\begin{equation}
\label{nt2}
 f(x)= f(x_n)+\frac{(x-x_n)}{2M}\sum_{k=1}^{2M} f'\big(x_n+\frac{(x-x_n)(k-0.5)}{2M}\big).
\end{equation}
Now, looking for $f(x)=0$ we arrive at
\begin{equation}
\label{nt3}
 x_{n+1}=x_n-\frac{2M(f(x_n))}{\sum_{k=1}^{2M} f'\big(x_n+\frac{(x-x_n)(k-0.5)}{2M}\big)}
\end{equation}
Further, substitute $x_{n+1}=x$ in (\ref{nm}) and replace $x-x_n$ in (\ref{nt3}) we obtain the new method as 
\begin{equation}\label{vnt}
 x_{n+1}=x_n-\frac{2M(f(x_n))}{\sum_{k=1}^{2M} f'\big(x_n- \frac{f(x_n)}{f'(x_n)}\frac{(k-0.5)}{2M}\big)}
\end{equation}
\noindent
\section{\normalsize\bf Convergence Analysis}
\begin{tdf}\label{T}
 Let the function $f : \D \subset \R \sra \R$ has a simple root $\alpha \in \D$, where $\D$ is an open interval. Assume
$f(x)$ has first, second and third derivatives in the interval $\D$. If the initial guess $x_0$ is closed to $\alpha$, then the method defined by (\ref{vnt}) converges cubically to
$\alpha$.
\end{tdf}
\noindent{\it Proof.} Let $\alpha$ is the simple root of $f(x)$ and $x_n= \alpha+e_n$. A use of Taylor expansion with $f(\alpha)=0$, we have
\begin{eqnarray}\label{nt3.6}
 f(x_n)=f'(\alpha)\big( e_n+C_2e_n^2+C_3e_n^3+O(e_n^4)\big),
\end{eqnarray}
where $C_k=\frac{1}{k!}\frac{f^k(\alpha)}{f'(\alpha)}$. Again
\begin{eqnarray}\label{nt3.7}
 f'(x_n)=f'(\alpha)\big(1+2C_2e_n+3C_3e_n^2+4C_4e_n^3+O(e_n^4)\big).
\end{eqnarray}
Further, dividing (\ref{nt3.6}) by (\ref{nt3.7}) yields
\begin{eqnarray}\label{nt3.8}
 \frac{f(x_n)}{f'(x_n)}=\big(e_n-C_2e_n^2+2(C_2^2-C_3)e_n^3+O(e_n^4)\big).
\end{eqnarray}
Now,
\begin{eqnarray}\label{nt3.9}
 x_n-M_k\frac{f(x_n)}{f'(x_n)}=x_n-M_k\big(e_n-C_2e_n^2+2(C_2^2-C_3)e_n^3+O(e_n^4)\big),
\end{eqnarray}  
where $M_k=\frac{(k-0.5)}{2M}$.\\
Equation (\ref{nt3.9}) rewrite as
\begin{align}\label{nt3.10}
 x_n-M_k\frac{f(x_n)}{f'(x_n)}&=x_n-M_k\big(e_n-C_2e_n^2+2(C_2^2-C_3)e_n^3+O(e_n^4)\big),\nonumber \\
 & = x_n-M_ke_n+M_kC_2e_n^2-2(C_2^2-C_3)M_ke_n^3+O(e_n^4),\nonumber\\
 &= x_n-(1-1+M_k)e_n+M_kC_2e_n^2-2(C_2^2-C_3)M_ke_n^3+O(e_n^4),\nonumber\\
 &= x_n-e_n+(1-M_k)e_n+M_kC_2e_n^2-2(C_2^2-C_3)M_ke_n^3+O(e_n^4),\nonumber\\
 &= \alpha+(1-M_k)e_n+M_kC_2e_n^2-2(C_2^2-C_3)M_ke_n^3+O(e_n^4).
\end{align}
From (\ref{nt3.10}) we can easily find that
\begin{align}\label{nt3.11}
 f'(x_n-M_k\frac{f(x_n)}{f'(x_n)})&=f'(\alpha)\big(1+2C_2(1-M_k)e_n+(2C_2^2M_k+3C_3(1-M_k)^2)e_n^2+O(e_n^3)\big).
\end{align}
Hence,
\begin{align}\label{nt3.12}
\sum_{k=1}^Nf'(x_n-M_k\frac{f(x_n)}{f'(x_n)})&=\sum_{k=1}^Nf'(\alpha)\big(1+2C_2(1-M_k)e_n+(2C_2^2M_k+3C_3(1-M_k)^2)e_n^2+O(e_n^3)\big),\nonumber\\
&=f'(\alpha)\sum_{k=1}^N\big(1+2C_2(1-M_k)e_n+(2C_2^2M_k+3C_3(1-M_k)^2)e_n^2+O(e_n^3)\big),\nonumber\\
&=f'(\alpha)\big(\sum_{k=1}^N1+2C_2e_n\sum_{k=1}^N(1-M_k)+2C_2^2e_n^2\sum_{k=1}^NM_k+3C_3e_n^2\sum_{k=1}^N(1-M_k)^2+O(e_n^3)\big),\nonumber\\
&=f'(\alpha)\big(N+2C_2Ne_n-2C_2e_nN/2+2C_2^2e_n^2N/2+3C_3e_n^2N\nonumber\\
&+3C_3e_n^2\big(\frac{N}{3}-\frac{1}{12N}\big)-6C_3e_n^2N/2+O(e_n^3)\big)\nonumber\\
&=f'(\alpha)\big(N+2C_2Ne_n+\big(NC_2^2+NC_3-\frac{C_3}{4N}\big)e_n^2+O(e_n^3)\big)
\end{align}
substitute (\ref{nt3.12}) and (\ref{nt3.6}) in(\ref{vnt}) we obtain
\begin{align}\label{nt3.13}
 x_{n+1}&=x_n-\frac{\big( e_n+C_2e_n^2+C_3e_n^3+O(e_n^4)\big)}{\big(N+2C_2Ne_n+\big(NC_2^2+NC_3-\frac{C_3}{4N}\big)e_n^2+O(e_n^3)\big)},\nonumber\\
 &=x_n- \big( e_n+\big(-C_2^2+\frac{C_3}{4N^2}\big)e_n^3+O(e_n^4)\big).
\end{align}
Substract $\alpha$ from both side of (\ref{nt3.13}), then we have
\begin{equation}
 e_{n+1}= \big(C_2^2-\frac{C_3}{4N^2}\big)e_n^3+O(e_n^4)
\end{equation}
this completes the rest of the proof. \hfill {$\Box$}

\section{\normalsize\bf Numerical Examples}
In this section, we present some numerical  results for various third order convergent iterative methods.
The following methods were compared:\\
Modified Newton's Method in Weerakoon and Fernando \cite{WF} $ x_{n+1}= x_n -\frac{2f(x_n)}{f'\big(xn-\frac{f(x_n)}{f'(x_n)}\big)+f'(x_n)}$.\\
\noindent
Modified Newton's Method proposed by Frontini et al. \cite{FS} $ x_{n+1}= x_n -\frac{f(x_n)}{f'\big(x_n-2\frac{f(x_n)}{f'(x_n)}\big)}$.\\
\noindent
Modified Newton's Method derived by Ozbal et al. \cite{OZ} $ x_{n+1}= x_n -\frac{f(x_n)}{2}\big(\frac{1}{f'(x_n)}+\frac{1}{f'\big(x_n-\frac{f(x_n)}{f'(x_n)}\big)}\big)$.\\
\noindent
Modified Newton's  Method derived by Kou et al. \cite{KLW} $ x_{n+1}= x_n -\frac{f\big(x_n+\frac{f(x_n)}{f'(x_n)}\big)-f(x_n)}{f'(x_n)}$.\\
\noindent
and method proposed by (\ref{vnt}) with $M=1$, $x_{n+1}=x_n-\frac{2(f(x_n))}{\sum_{k=1}^{2} f'\big(x_n- \frac{f(x_n)}{f'(x_n)}\frac{(k-0.5)}{2}\big)}$.\\
For every problem an attempt made to find an approximation $x_n$ of the simple root of equation $f(x) = 0$ through $n$ times
iterations. The number of function evaluations (NFE) is counted as the sum of the number of evaluations of the function $f$ and its 
first order derivative $f'$ . The computational results are displayed in Table 1.\\
The numerical experiments carried over the following equation:
\begin{eqnarray}
 f_1(x)&=&x^5-x+1,\nonumber\\
 f_2(x)&=&cosx-x ,\nonumber\\
 f_3(x)&=&arctanx,\nonumber\\
 f_4(x)&=&10xe^{-x^2}-1,\nonumber\\
 f_5(x)&=&e^{-x}sinx+log(x^2+1),\nonumber\\
 f_6(x)&=&x^3-e{-x},\nonumber\\
 f_7(x)&=&e{-x}-cosx,\nonumber\\
\end{eqnarray}
\noindent
The numerical results presented in Table 1 show that the proposed method has perform equally as compared with the 
other methods of the same order. Thus, the new methods can compete with other third-order methods in literature.\\
\begin{table}[!htb]
\centering
\begin{tabular}{|l|l|l|l|l|l|}
\hline
Function & $ x_0$&Various Method & IT & NFE  & $x_n$     \\
    &       &    &         &     &                  \\
\hline
\hline $f_1$     &  2  &  MNM(\cite{WF})   &  &       & Diverse \\
                 &     &  MNM(\cite{FS})   &  &        & Diverse \\
                 &     &  MNM(\cite{OZ})   & 13&   39     & -1.16730397826142\\
                 &     &  MNM(\cite{KLW})  & 17 &  51      & -1.16730397826142\\
                 &     &  MNM New          & 9 &  36     &-1.16730397826142\\                 
\hline $f_2 $    &  1.2& MNM(\cite{WF})  & 4 & 12 & 0.739085133215161  \\
                 &     & MNM(\cite{FS})  &  &        & Diverse\\
                 &     &  MNM(\cite{OZ})    &  &        & Diverse\\
                 &     &  MNM(\cite{KLW})     &  4& 12       & 0.739085133215161\\
                 &     &  MNM New      &  4& 16       & 0.739085133215161\\ 
\hline$ f_3$     &  3  &  MNM(\cite{WF}) &  &  & Diverse  \\
                 &     & MNM(\cite{FS}) &  &        & Diverse\\
                 &     &  MNM(\cite{OZ})    &  &        & Diverse\\
                 &     &  MNM(\cite{KLW})     & 4 & 12       & 0.000000000000015\\
                 &     &  MNM New      & 4 &  16      & 0.0\\ 
\hline$f_4 $     & 2.5 & MNM(\cite{WF})   & 4 & 12 &1.67963061042845   \\
                 &     &  MNM(\cite{FS})  & 7 & 21 &1.67963061042845 \\
                 &     & MNM(\cite{OZ})   & 4 & 12 & 0.101025848315685\\
                 &     & MNM(\cite{KLW})  & 6 & 18 & 1.67963061042845\\
                 &     &  MNM New         & 5 & 20 & 1.67963061042845\\ 
\hline$ f_5$ & 1.3 & MNM(\cite{WF})  & 3& 9   & 7.68481808334733E-021  \\
                 &                   &  MNM(\cite{FS})  & 4 & 12 &1.7197167733818E-028\\
                 &                   & MNM(\cite{OZ})    & 4&  12 & 5.12759588393657E-030\\
                 &                   & MNM(\cite{KLW})      &  3&  9 & 3.82180552357862E-019\\
                 &                   &  MNM New      & 3 & 12 & 4.53468286561001E-017\\ 
\hline$ f_6$ & 2 &   MNM(\cite{WF})  & 5 &  15     & 0.77288295914921        \\
                 &     &  MNM(\cite{FS})   & 4 &  12     & 0.77288295914921\\
                 &     &  MNM(\cite{OZ})     &  4&   12     & 0.77288295914921\\
                 &     &  MNM(\cite{KLW})       &  5&   15     & 0.77288295914921\\
                 &     &  MNM New       &  4&  16     & 0.77288295914921\\ 
\hline$ f_7$ &2 &  MNM(\cite{WF})  & 3  &  9     & 1.29269571937339         \\
                 &     &  MNM(\cite{FS})   & 3 &  9     &1.29269571937339 \\
                 &     &  MNM(\cite{OZ})     & 4 &  12    & 1.29269571937339\\
                 &     & MNM(\cite{KLW})       & 4 &  12    & 1.29269571937339\\
                 &     &  MNM New       & 3 &    12   & 1.29269571937339\\ 
\hline

\end{tabular}
\caption{Comparison of various third order convergent iterative methods and the New Newton method}
\end{table}
\section{\normalsize\bf Conclusion}
This article deals with a new modified Newton methods for solving nonlinear equations.
In Theorem (\ref{T}), it is proved that the new method has third order convergence. Further, the new methods can compete with other
third-order methods. Finally, numerical experiments are conducted to confirm our theoretical findings.
The new method has great practical utility.

\end{document}